\documentstyle[amssymb, amsmath]{article}
\begin{document}
\newcommand{\bref}[1]{$\mbox{(\ref{#1})}$}
\newcommand{\w}{\sqrt {t + b}}
\newcommand{\p}{{\Bbb {R^+}} \rightarrow {\Bbb {C}}}
\newcommand{\n}{\sqrt {2n + 1}}
\newcommand{\r}{{\Bbb {R}} ^d}
\newcommand{\q}{{\mathcal Q}}
\newcommand{\z}{{\Bbb {Z}} ^d}
\newtheorem{thm}{Theorem}
\newtheorem{lem}{Lemma}
\newtheorem{rem}{Remark}
\newtheorem{cor}{Corollary}

\title 
{{\bf Estimates of functions with vanishing periodizations
}}

\author{Oleg Kovrizhkin \thanks{This research was partially conducted by the 
author for the Clay Mathematics Institute and partially supported by NSF grant 
DMS 97-29992}}

\date{}

\maketitle

\begin{center}
\rm 
MIT, 2-273, Dept of Math\\
77 Mass Ave\\
Cambridge, MA 02139, USA\\
{\sl E-mail address}: oleg@math.mit.edu\\
\end{center}

\begin{abstract} We prove that if a function $f \in L^p({\Bbb {R}} ^d)$ has 
vanishing periodizations then $\|f\|_{p'} \lesssim \|f\|_{p}$, provided $1 \le p 
< \frac {2d}{d + 2}$ and dimension $d \ge 3$.
\end{abstract}

\section{Introduction}\label{intro}

Let $f \in L^1({\Bbb {R}} ^d)$. Define a family of its periodizations with 
respect to a rotated integer lattice:
\begin{eqnarray} g_{\rho}(x) = \sum \limits _{\nu \in {\Bbb {Z}} ^d} f(\rho (x - 
\nu))\label{g}\end{eqnarray}
for all rotations $\rho \in SO(d)$. We have a trivial estimate $\|g_{\rho}\|_1 
\le \|f\|_1$ and $\widehat{g_{\rho}}(m) = \hat f (\rho m)$ where $m = (m_1, ..., 
m_d) \in \z$. The author has shown recently that $g_{\rho}$ is in $L^2([0,1]^d 
\times SO(d))$ if and only if $f \in L^2(\r)$, provided the dimension $d \ge 5$. 
The requirement $f \in L^1({\Bbb {R}} ^d)$ can be replaced by $f \in L^p({\Bbb 
{R}} ^d)$ for a certain range of $p$, see for details (\cite{OK1}), 
(\cite{OK2}).\\

The main object of our research will be functions $f$ whose periodizations 
$g_{\rho}$ identically vanish for a.e. rotations $\rho \in SO(d)$. It is 
equivalent to the statement that $\hat f$ vanishes on all spheres of radius $|m| 
= (m_1^2 + ... + m_d^2)^{\frac 12}$ where $m \in \z$. Such functions are closely 
related to the Steinhaus tiling set problem (\cite{K}),
(\cite{KW}): does there exists a (measurable) set $E \subset \r$ such that every
rotation and translation of $E$ contains exactly one integer lattice point? 
M. Kolountzakis (\cite{K}) showed that if $f \in L^1$ and 
$|x|^{\alpha}f(x) \in L^1$ for a certain $\alpha > 0$ and $f$ has constant 
periodizations then $\hat f \in L^1$ when dimension $d = 2$. M. Kolountzakis and 
T. Wolff (\cite{KW}, Theorem 1) proved that if periodizations of a function from 
$L^1({\Bbb {R}} ^d)$ are constants then the function is continuous and, in fact, 
bounded, provided that the dimension d is at least three.  We will generalize 
the 
last result for functions $f$ in $L^1({\Bbb {R}} ^d) \cap L^p({\Bbb {R}} ^d)$:
\begin{thm}
Let $d \ge 3$ and $f \in L^1({\Bbb {R}} ^d) \cap L^p({\Bbb {R}} ^d)$, $1 \le p < 
\frac {2d}{d + 2}$, has identically vanishing periodizations then $f \in 
L^{p'}({\Bbb {R}} ^d)$:
\begin{eqnarray*} \|f\|_{p'} \le C\|f\|_p \end{eqnarray*}
where $C$ depends only on $d$ and $p$.
\end{thm}
The main reason why the dimension $d \ge 3$ comes from the famous Lagrange theorem saying that every positive integer can be represented as sums of four 
squares and actually from the fact that every integer of form $8k + 1$ can be written as sums of three squares. Since relatively few integers can be 
represented as sums of two squares, we will show in {\bf Section 3} that the result 
of M. Kolountzakis and T. Wolff doesn't hold if $d = 2$ and that is why there is no theorem for $d = 2$. Another reason why the dimension $d \ge 3$ is because 
we consider the family of periodizations with respect to the $SO(d)$ group of rotations. It leads to estimates involving the decay of spherical harmonics. The 
rate of decay for $d = 2$ is not fast enough although it is almost fast enough. That is why for $d = 2$ the range of $p$ in the theorem becomes empty: $1 \le p < 1$. 
\begin{rem}
There is no essential difference between the case of identically vanishing 
periodizations and the case of $g_{\rho}$ being trigonometric polynomials of 
uniformly bounded degrees for all $\rho \in SO(d)$.
\end{rem}
\begin{cor}
If $p \le r \le p'$
then under the conditions of {\bf Theorem 1}
\begin{eqnarray*} \|f\|_{r} \le C\|f\|_p \end{eqnarray*}
where $C$ depends only on $d$ and $p$.
\end{cor}
We will show in {\bf Section 3} that this range of $r$ is sharp.\\

We will use the notation $x \lesssim y$ meaning 
$x \le Cy$, and $x \sim y$ meaning that $x \lesssim y$ and $y \lesssim x$ for 
some constant 
$C > 0$ independent from $x$ and $y$.\\

\section{Proof of the theorem}\label{theor}
Define the following functions $h, h_1, h_2: {\r \times \Bbb {R^+}} \rightarrow 
{\Bbb {C}}$
\begin{eqnarray}h(y, t) &=& \int \hat f(\xi)e^{i2\pi y\cdot\xi} 
d\sigma_t(\xi)\label{h}\\
&=& \int_{\r} f (x) \widehat {d\sigma_t}(y - x)dx\\
&=& \int_{\r} f (y - x) \widehat {d\sigma_t}(x)dx, \label{h3}\end{eqnarray}
\begin{eqnarray}h_1(y, t) 
= \int_{|x| \le 1} f(y - x) \widehat {d\sigma_t}(x)dx, \label{h1}\end{eqnarray}
\begin{eqnarray}h_2(y, t) 
= \int_{|x| > 1} f (y - x) \widehat {d\sigma_t}(x)dx \label{h2}\end{eqnarray}
where $d\sigma_t$ is the Lebesgue surface measure on a sphere of radius $t$.
Clearly, $h = h_1 + h_2$. To proceed further we will need certain technical estimates associated with $h_1$ and $h_2$ proven in two lemmas below. The proof of 
the theorem itself starts after {\bf Remark 2} to {\bf Lemma 2}. The Fourier transforms in these two lemmas below are taken 
with respect to variable $t$, except in the second part of the proof of {\bf 
Lemma 2}. $L^{p'}$ norms are taken over variable $y$. We will  apply some 
technique from M. Kolountzakis and T. Wolff (\cite{KW}) and O. 
Kovrijkine (\cite {OK1}, \cite {OK2}).

~

\begin{lem} 
Let $q: {\Bbb {R}} \rightarrow {\Bbb {R}}$ be a Schwartz function supported in 
$[\frac 12, 2]$, let $f \in L^p(\r)$ where $1 \le p \le 2$ and let $b \in 
[0,1)$. Define $H_{1, N}:{\r \times \Bbb {R}} \rightarrow {\Bbb {C}}$
\begin{eqnarray*} H_{1, N} (y, t) = \frac 1 {\w} h_1(y, \w) q(\frac {\w} 
N).\end{eqnarray*}
Then  
\begin{eqnarray} \sum \limits _{l \ge 0}\sum \limits _{\nu \neq 0}\|\hat H_{1, 
2^l} (y, \nu) \|_{p'} \le C\|f\|_p \label{H1}\end{eqnarray}
where $C$ depends only on $q$ and $d$.
\end{lem}
{\bf Proof of Lemma 1:}\\
It will be enough to show that
\begin{eqnarray} \sum \limits _{\nu \neq 0}\|\hat H_{1, N} (y, \nu) \|_{p'} \le 
\frac {C\|f\|_p} {N} . \label{H11}\end{eqnarray}
We have
\begin{eqnarray} |\hat H_{1, N}(y, \nu)| &\le& \frac {C} {|\nu| ^k}\int|\frac 
{\partial ^ k}{\partial t^k}H_{1, N}(y, t)|dt \label{Premin}
\end{eqnarray}
for $\nu \neq 0$.
Applying Minkowski's inequality to \bref{Premin} we have
\begin{eqnarray}
\|\hat H_{1, N} (y, \nu) \|_{p'} &\le& 
\frac {C} {|\nu| ^k}\int\|\frac {\partial ^ k}{\partial t^k}H_{1, N}(y, 
t)\|_{L^{p'}(dy)}dt. \label{Aftermin}
\end{eqnarray}
We need to estimate the integrand on the right side of \bref{Aftermin}.
To do so we will first estimate the $L^{p'}$ norm of derivatives of $h_1(y, t)$ 
when
$t \ge 1$:
\begin{eqnarray}\|\frac {\partial ^ k}{\partial t^k}h_1(y, t)\|_{p'} \lesssim 
t^{d-1}\|f\|_p\label{dif}\end{eqnarray}
with an implicit constant depending  only on $k$ and $d$.
In order to obtain  
\bref{dif}, rewrite the definition of $h_1$ \bref{h1} in the following way:
\begin{eqnarray*}h_1 (y, t) &=& \int_{|x| \le 1} f(y - x) \widehat 
{d\sigma_t}(x)dx\\
&=&t^{d-1}\int_{\r} f (y - x)\cdot\chi_{\{|x| \le 1\}}\int _{|\xi| = 1} 
e^{-i2\pi t x\cdot \xi}d\sigma(\xi)dx, \end{eqnarray*}
 differentiate the last equality $k$ times and apply Young's inequality.\\
 
We can easily prove by induction that
\begin{eqnarray}\frac {d^k} {dt^k} \left (\frac {h_1 (\w)} {\w} \right) =
 \sum \limits _{i = 0} ^ {k} C_{i,k}\frac {h_1^{(i)}(\w)}{(\w)^{2k + 1 
-i}}.\label {dif1}\end{eqnarray}
Combining \bref{dif1} and \bref{dif} we obtain for $t \thicksim N^2$
\begin{eqnarray}\left \| \frac {\partial ^k} {\partial t^k} \left (\frac {h_1 
(y, \w)} {\w} \right) \right \|_{p'}
 \le CN^{d-k-2}\|f\|_p\label{dif3}\end{eqnarray}
with $C$ depending only on $k$ and $d$.\\

Since $q(\frac {(\w)} {N}) = q(\sqrt {t' + b'}) = \tilde q(t')$ with $t' = \frac 
t {N^2}$
 and $b' = \frac b {N^2}$ and $\tilde q(t')$ is a Schwartz function supported in
 $t' \thicksim 1$, we have 
\begin{eqnarray}\left |\frac {d^k} {dt^k} q(\frac {(\w)} {N}) \right| &=&
N^{-2k}|\frac {d^k} {d{t'}^k} \tilde q (t')|\nonumber \\
&\le& CN^{-2k}\label {dif2}\end{eqnarray}
with $C$ depending only on $k$ and $q$.\\

$q(\frac {(\w)} {N})$ is supported in $t \thicksim N^2$ hence we obtain from 
\bref{dif3} and \bref{dif2} that
\begin{eqnarray}\left \|\frac {\partial ^k} {\partial t^k} H_{1, N}(y, t)\right 
\|_{p'} &=& 
\left \|\frac {d^k} {dt^k} \left (\frac {h_1 (y, \w)} {\w} q(\frac {\w} 
{N})\right) \right\|_{p'}
\nonumber \\
 &\le& CN^{d-2-k}\|f\|_p\label{dif4}\end{eqnarray}
with $C$ depending only on $k$, $d$ and $q$.  Since $H_{1, N}(y, t)$ is also 
supported in 
$t \thicksim N^2$ we have
$$\int\|\frac {\partial ^ k}{\partial t^k}H_{1, N}(y, t)\|_{L^{p'}(dy)}dt \le 
CN^{d-k}\|f\|_p.$$
Substituting the above estimate to \bref{Aftermin} we obtain
\begin{eqnarray}
\|\hat H_{1, N} (y, \nu) \|_{p'} &\le& 
\frac {CN^{d-k}\|f\|_p} {|\nu| ^k}\label{nu}\end{eqnarray}
for every $\nu \neq 0$.\\

Summing \bref{nu} over all $\nu \neq 0$ and putting $k = d + 1$ we get our 
desired result
\begin{eqnarray*} \sum \limits _{\nu \neq 0}\|\hat H_{1, N} (y, \nu) \|_{p'} \le 
\frac {C\|f\|_p} {N} . \end{eqnarray*}
where $C$ depends only on $q$ and $d$. Sum over dyadic $N$ to obtain the 
statement of the lemma. \hfill$\square$

~

The next lemma will be proven in the spirit of the Stein-Tomas restriction theorem (\cite 
{DC}, p.104).
\begin{lem}
Let $q: {\Bbb {R}} \rightarrow {\Bbb {R}}$ be a Schwartz function supported in 
$[\frac 12, 2]$, let $f \in L^p({\Bbb {R}} ^d)$ where $1 \le p < \frac {2d} {d + 
2}$ and let $b \in [0,1)$. Define $H_{2,N}:{\r \times \Bbb {R}} \rightarrow 
{\Bbb {C}}$
\begin{eqnarray*} H_{2, N} (y, t) = \frac 1 {\w} h_2(y, \w) q(\frac {\w} 
N).\end{eqnarray*}
Then we have 
\begin{eqnarray} \sum \limits _{\nu \neq 0}\|\sum \limits _{l \ge 0} \hat H_{2, 
2^l} (y, \nu)\|_{p'} \le C\|f\|_p\label{H2}\end{eqnarray}
with $C$ depending only on $p$, $q$ and $d$.
\end{lem}
{\bf Proof of Lemma 2:}\\
We have
\begin{eqnarray}&&\hat H_{2, N} (y, \nu)\nonumber \\
 &=& \int H_{2, N}(y, t)e^{-i2\pi \nu t}dt\nonumber \\
 &=& 2e^{i2\pi \nu b}\int Nq(t) h_2(y, tN)e^{-i2\pi \nu (Nt)^2}dt\nonumber \\
&=& 2e^{i2\pi \nu b}\int Nq(t)e^{-i2\pi \nu (Nt)^2}\int_{|x| > 1} f (y - x) 
\widehat {d\sigma_{Nt}}(x)dxdt\nonumber \\
&=& 2e^{i2\pi \nu b}\int_{|x| > 1} f (y - x)\int Nq(t)e^{-i2\pi \nu (Nt)^2} 
(Nt)^{d-1}\widehat {d\sigma}(Ntx)dtdx \label{H3} \\
&=& (D_{N, \nu} * f) (y)
\nonumber \end{eqnarray}
where 
\begin{eqnarray} D_{N, \nu}(x) = 2e^{i2\pi \nu b}\chi_{\{|x| > 1\}}\int 
Nq(t)
e^{-i2\pi \nu (Nt)^2} (Nt)^{d-1}\widehat 
{d\sigma}(Ntx)dt.\label{D99}\end{eqnarray}
Denote by 
\begin{eqnarray} K_{\nu}(x) = \sum \limits _{l \ge 0}D_{2^l, \nu}(x) 
\label{K99}.
\end{eqnarray}
We need to estimate
\begin{eqnarray*}\|\sum \limits _{l \ge 0} \hat H_{2, 
2^l} (y, \nu)\|_{p'} = \|K_{\nu} * f\|_{p'}.
\end{eqnarray*}
If $p' = \infty$ or $p' = 2$ we have
\begin{eqnarray*} \|K_{\nu}*f\|_{\infty} \le \|K_{\nu}\|_{\infty}\|f\|_1 \\
\|K_{\nu}*f\|_{2} \le \|\hat K_{\nu}\|_{\infty}\|f\|_2.\end{eqnarray*}
First we will show that 
\begin{eqnarray} \|K_{\nu}\|_{\infty} &\le& \|\sum \limits _{l \ge 0}
|D_{2^l, \nu}|(x)\|_{\infty} \nonumber \\
&\le& C|\nu| ^{-\frac{d}2}\label{K1}.\end{eqnarray}
To do so we need to estimate $D_{N, \nu}$.

We will use a well-known fact that $\widehat {d\sigma} (x) = Re(B(|x|))$ with 
$B(r)
 = a(r)e^{i2\pi r}$ and $a(r)$ satisfying estimates
\begin{eqnarray}|a^{k}(r)| \le \frac {C}{r^{\frac{d-1}2 + 
k}}\label{a}\end{eqnarray}
with $C$ depending only on $k$ and $d$.
Now we will  estimate the  integral in \bref{D99} with $B(|x|)$ 
instead of
 $\widehat {d\sigma}(x)$
\begin{eqnarray}&&\int Nq(t)e^{-i2\pi \nu (Nt)^2}(Nt)^{d-1} a(N|x|t)e^{i2\pi 
N|x|t }dt
\nonumber \\ &=&\frac {N^{\frac {d+1} 2}} {|x|^{\frac {d-1} 2}}\int q(t)
e^{-i2\pi \nu (Nt)^2}t^{d-1} a(N|x|t)(N|x|)^{\frac {d-1} 2}e^{i2\pi N|x|t 
}dt\nonumber \\
&=&\frac {N^{\frac {d+1} 2}} {|x|^{\frac {d-1} 2}}e^{i2\pi \frac 
{|x|^2}{4\nu}}\int q(t) a(N|x|t)(N|x|)^{\frac {d-1} 2}t^{d-1}e^{-i2\pi \nu N^2 
(t-\frac {|x|}{2\nu N})^2}dt\nonumber \\
&=&\frac {N^{\frac {d+1} 2}} {|x|^{\frac {d-1} 2}}e^{i2\pi \frac 
{|x|^2}{4\nu}}\int \phi(t, |x|)e^{-i2\pi \nu N^2 (t-\frac {|x|}{2\nu N})^2}dt
\label{a1}\end{eqnarray}
where $\phi(t, |x|) = q(t) a(N|x|t)(N|x|)^{\frac {d-1} 2}t^{d-1}$ is a Schwartz 
function
with respect to variable $t$ supported in $[\frac 12, 2]$ which is bounded, together with each derivative
uniformly in $t$, $|x| \ge 1$ and $N$ because of \bref{a}. Note that we used 
here the fact
that $N|x| \ge 1$. We can say even more. 
 Let $|x| = c\cdot r$ where $c \ge 2$ and $r \ge \frac 1 2$. Then all partial 
derivatives
 of $\phi(t, c\cdot r)$ with respect to $t$ and $r$ are also bounded uniformly 
in $t$,
$r$, $c$ and $N$. Hence $\phi(t, c\cdot t)$ is a Schwartz function supported in 
$[\frac 12, 2]$ which is bounded, together with each derivative
uniformly in $t$, $c$ and $N$. We will use this fact later to estimate $\hat 
K_{\nu}$.\\

 Fix some $|x| \ge 1$. In the calculations below we will write just $\phi(t)$ 
instead of $\phi(t, |x|)$ for simplicity.
From the method of stationary phase (\cite{H}, Theorem 7.7.3) it follows that if 
$k \ge 1$ then
\begin{eqnarray}|\int \phi(t)e^{-i2\pi \nu N^2 (t-\frac {|x|}{2\nu N})^2}dt - 
\sum \limits _{j=0}^{k-1}c_j(\nu N^2)^{-j-\frac 12}\phi^{(2j)}(\frac {|x|}{2\nu 
N})|
 \le c_k(|\nu| N^2)^{-k-\frac 12}\label{phi}\end{eqnarray}
where $c_j$ are some constants.

Since $\phi$ is supported in $[\frac 12, 2]$ we conclude from \bref{phi} that 
\begin{eqnarray}|\int \phi(t)e^{-i2\pi \nu N^2 (t-\frac {|x|}{2\nu N})^2}dt| \le 
\begin{cases} C(|\nu| N^2)^{-\frac 12}&\text{if $N \in [\frac {|x|} {4\nu},\frac 
{|x|} {\nu}]$}\\
C_k(|\nu| N^2)^{-k-\frac 12}&\text{if $N \notin [\frac {|x|} {4\nu},\frac {|x|} 
{\nu}]$}\end{cases}
.\label{phi1}\end{eqnarray}
Replacing in \bref{D99} $\widehat {d\sigma}(x)$ with $\frac {B(|x|) + \bar 
B(|x|)} 2$ it follows from \bref{phi1} that
\begin{eqnarray} |D_{N, \nu}(x)| \le \frac{N^{\frac{d+1} 2}}{|x|^{\frac{d-1} 
2}}\begin{cases} C(|\nu| N^2)^{-\frac 12}&\text{if $N \in [\frac {|x|} 
{4|\nu|},\frac {|x|} {|\nu|}]$}\\
C_k(|\nu| N^2)^{-k-\frac 12}&\text{if $N \notin [\frac {|x|} {4|\nu|},\frac 
{|x|} 
{|\nu|}]$}
\end{cases}
.\label{D1}\end{eqnarray}
The number of dyadic $N \in [\frac {|x|} {4\nu},\frac {|x|} 
{\nu}]$ is
 at most $3$. Therefore choosing $k \ge \frac {d-1}2$ and summing \bref{D1} 
over all dyadic $N$ we have
\begin{eqnarray*}|K_{\nu}(x)| \le \sum \limits _{l \ge 0}|D_{2^l, \nu}(x)| \le 
C|\nu|^{-\frac {d} 
2}
\end{eqnarray*}
with $C$ depending only on $d$ and $q$. Thus we proved \bref{K1}.\\
Now we will show that
\begin{eqnarray} \|\hat K_{\nu}\|_{\infty} \le \|\sum \limits _{l \ge 0}|\hat 
D_{2^l, \nu}|(y)\|_{\infty} \le C \label{K2}.\end{eqnarray}
Since supp $\phi \in [\frac 12, 2]$ we can re-write \bref{phi} for a stronger 
version of 
the method of stationary phase (\cite{H}, Theorems 7.6.4, 7.6.5, 7.7.3)
\begin{eqnarray}|\int \phi(t)e^{-i2\pi \nu N^2 (t-\frac {|x|}{2\nu N})^2}dt - 
\sum \limits _{j=0}^{k-1}c_j(\nu N^2)^{-j-\frac 12}\phi^{(2j)}(\frac {|x|}{2\nu 
N})| \le \frac 
{c_k(|\nu| N^2)^{-k-\frac 12}} {max(1, \frac {|x|} {8N|\nu|})^k} 
\nonumber\end{eqnarray}
where $c_j$ are some constants.
Therefore, if $\nu > 0$, 
\begin{eqnarray}D_{N, \nu}(x) = \chi_{\{|x| > 1\}}\frac {N^{\frac {d+1} 2}} 
{|x|^{\frac {d-1} 2}}e^{i2\pi \frac {|x|^2}{4\nu}}\sum \limits 
_{j=0}^{k-1}c_j(\nu N^2)^{-j-\frac 12}\phi^{(2j)}(\frac {|x|}{2\nu N}) + 
\phi_k(x)\label{D10}\end{eqnarray}
where $|\phi_k(x)| \le \chi_{\{|x| > 1\}}\frac {N^{\frac {d+1} 2}} {|x|^{\frac 
{d-1} 2}}\frac {c_k(|\nu| N^2)^{-k-\frac 12}} {max(1, \frac {|x|} 
{8N|\nu|})^k}$.
If $\nu < 0$ then just replace $\phi^{(2j)}(\frac {|x|}{2\nu N})$ with 
$\bar{\phi}^{(2j)}(-\frac {|x|}{2\nu N})$. We further assume that $\nu > 0$.
Choosing $k \ge \frac {d + 2} 2$ we have 
\begin{eqnarray}\|\hat \phi_k\|_{\infty} &\le& \|\phi_k\|_1\nonumber \\
&=& \int\limits_{|x| \le 8\nu N}|\phi_k| dx + \int\limits_{|x| > 8\nu N}|\phi_k| 
dx \nonumber \\
&\le& \frac C N \label{phik}\end{eqnarray}
where $C$ depends only on $d$ and $q$.  We can ignore $\chi_{\{|x| > 1\}}$ in 
front of the sum in 
\bref{D10} because if $\frac {|x|}{2\nu N} \in [\frac 12, 2]$, then $|x| \ge \nu 
N \ge 1$.
 We will consider only the zero term in the sum. The other terms can be treated 
similarly. 
The Fourier transform of 
\begin{eqnarray*}\frac {N^{\frac {d+1} 2}} {|x|^{\frac {d-1} 2}}e^{i2\pi \frac 
{|x|^2}{4\nu}}
(\nu N^2)^{-\frac 12}\phi(\frac {|x|}{2\nu N})\end{eqnarray*}
at point $y$ is equal to
\begin{eqnarray}N^{\frac {d+1} 2}(2\nu N)^{\frac {d+1} 2}(\nu N^2)^{-\frac 
12}\int_{\r} \psi(|x|)e^{i2\pi \nu N^2 |x|^2} e^{-i2\pi 2\nu N x \cdot y} dx = 
\nonumber \\
C(\nu N^2)^{\frac d 2}e^{-i2\pi \nu|y|^2}\int_{\r} \psi(|x|)e^{i2\pi \nu N^2 |x 
- \frac y N|^2}dx\label{psi10}\end{eqnarray}
where $\psi(t) = \phi(t, 2\nu N t) t^{-\frac{d-1}2}$  is a Schwartz function 
supported in 
$[\frac 12, 2]$ whose derivatives and the function itself are bounded uniformly 
in $t$, 
$\nu$ and $N$ (see remark after \bref{a1}). The same is true about partial 
derivatives of $\psi(|x|)$. Applying the stationary phase method for $\r$ 
(\cite{H}, Theorem 7.7.3) we get 
\begin{eqnarray}|\int_{\r} \psi(|x|)e^{i2\pi \nu N^2 |x - \frac y N|^2}dx| \le 
\begin{cases} C(\nu N^2)^{-\frac d 2}&\text{if $N \in [\frac {|y|} {2},2|y|]$}\\
C_k(\nu N^2)^{-k-\frac d 2}&\text{if $N \notin [\frac {|y|} 
{2},2|y|]$}\end{cases}
.\label{psi1}\end{eqnarray}
Therefore the absolute value of \bref{psi10} can be bounded from above by:
\begin{eqnarray} \le \begin{cases} C&\text{if $N \in [\frac {|y|} {2},2|y|]$}\\
C_k(\nu N^2)^{-k}&\text{if $N \notin [\frac {|y|} {2},2|y|]$}\end{cases}
.\label{psi2}\end{eqnarray}
Similar inequalities hold for Fourier transforms for the rest of the terms in 
the sum in \bref{D10}.
The number of dyadic $N \in [\frac {|y|} {2},2|y|]$ is bounded by $3$. Using 
\bref{phik},
 choosing $k \ge 1$ in \bref{psi2} and summing over all dyadic $N$ we get
\begin{eqnarray}\sum \limits_{l \ge 0}|\hat D_{2^l, \nu}(y)| \le C 
\label{D5}\end{eqnarray}
with $C$ depending only on $d$ and $q$, provided $\nu \neq 0$.
Thus we proved 
\bref{K2}.\\
Using \bref{K1} and \bref{K2} and interpolating between $p=1$ and $p=2$, we 
obtain
\begin{eqnarray} \|K_{\nu}*f\|_{p'} \le C|\nu| ^{-\alpha 
_p}\|f\|_p \label{H20}\end{eqnarray}
where $\alpha_p = \frac d 2 \frac {2-p}p$. $\alpha_p > 1$ if $p < \frac {2d}{d + 
2}$. 
Summing \bref{H20} over all $\nu \neq 0$, we get the desired inequality
\begin{eqnarray*} \sum \limits _{\nu \neq 0}\|\sum \limits _{l \ge 0} \hat H_{2, 
2^l} (y, \nu)\|_{p'} \le C\|f\|_p.\end{eqnarray*}\hfill$\square$

\begin{rem}
It is clear from the proof that we have the same inequality if the summation 
over $l 
\ge 0$ is replaced by summation over any subset of nonnegative integers.
\end{rem}

Now we are in a position to proceed with the proof of the theorem.
Let $q: {\Bbb {R}} \rightarrow {\Bbb {R}}$ be a fixed nonnegative Schwartz 
function 
supported in $[\frac 12, 2]$ such that 
$$q(t) + q(t/2) = 1$$
when $t \in [1,2]$. It follows that
\begin{eqnarray}\sum \limits_{l \ge 0} q(\frac t {2^l}) = 1 
\label{q10}\end{eqnarray}
when $t \ge 1$.
Define $$q_0(t) = 1 - \sum\limits_{l \ge 0}q(\frac t {2^l})$$ for $t \ge 0$. It 
is 
clear that $q_0(|x|)$ is a Schwartz function supported in $|x| \le 1$. Let 
$\psi(t) = q_0(t) + q(t)$ then 
$$\psi_k (t) = \psi(\frac t {2^k}) = q_0(t) + \sum\limits_{l 
\ge 0}^{k}q(\frac t {2^l})$$ and $\psi(|x|)$ is a Schwartz function supported in 
$|x| \le 2$ such that $\psi(|x|) = 1$ if $|x| \le 1$. Therefore
\begin{eqnarray*}\int \hat f (x) e^{2\pi x \cdot y}\psi(\frac {|x|} {2^k})dx = 
(f*\widehat {\psi_k})(y)\end{eqnarray*}
converges to $f$ in $L^p$ as $k \rightarrow \infty$. To prove that $f \in 
L^{p'}$ and $\|f\|_{p'} \lesssim \|f\|_p $ it will be enough to show that
\begin{eqnarray*}\|f*\widehat {\psi_k} 
\|_{p'} \le C\|f\|_p \end{eqnarray*}
since the claim will follow by an application of Fatou's lemma to a subsequence
of $f*\widehat {\psi_k}$ converging a.e. to $f$.\\

We have 
\begin{eqnarray}(f*\widehat {\psi_k})(y) 
&=& (f*\widehat{q_0})(y) +  \sum\limits_{l 
\ge 0}^{k}\int \hat f (x) e^{2\pi x \cdot y}q(\frac {|x|} {2^l})dx\nonumber \\
&=& (f*\widehat{q_0})(y) + \sum\limits_{l 
\ge 0}^{k}\int\limits_{0}^{\infty} q(\frac {t} {2^l}) \int \hat f(\xi)e^{i2\pi 
y\cdot\xi} 
d\sigma_t(\xi)dt\nonumber \\
&=& (f*\widehat{q_0})(y) + \sum\limits_{l 
\ge 0}^{k}\int\limits_{0}^{\infty} q(\frac {t} {2^l}) h(y,t) dt 
\label{convol}.\end{eqnarray}
Applying Young's inequality we estimate the first term:
\begin{eqnarray}\|f*\widehat{q_0}\|_{p'} \lesssim \|f\|_p 
\label{zeroter}\end{eqnarray}
for $1 \le p \le 2$. Now we have to estimate the sum over $l$. 

It is a well-known fact from Number Theory proven by Lagrange that every positive integer can be 
repersented as sums of four squares (\cite {G}, p.25), moreover there exists an infinite 
arithmetic progression of positive integers, e.g., $8n + 1$, which can be 
represented as sums of three squares (\cite {G}, p. 38). We will use only the latter fact. Therefore, rescaling we can 
assume that 
$\hat f$ vanishes on all spheres of radius $\sqrt{n + b}$ where $n$ is a 
nonnegative integer and $0 < b < 1$ is a 
fixed number. Therefore $h(y, \sqrt{n + b}) = 0$ for all $y \in \r$.
Making a change of variables and keeping in mind that $q$ is supported in 
$[\frac 1 2, 2]$ we re-write every term in the sum in the following way:
\begin{eqnarray*}\int\limits_{0}^{\infty} q(\frac {t} {N}) h(y,t) dt =
\int \frac 1 {2\w} q(\frac {\w} {N}) h(y,\w) dt.\end{eqnarray*}
 An application of Poisson's summation formula gives us
 \begin{eqnarray*}0 &=& \sum\limits_{n} \frac 1 {\sqrt{n + b}} q(\frac {\sqrt{n 
+ b}} {N}) h(y,\sqrt{n + b})\\
 &=& \sum\limits_{\nu} \left (\frac 1 {\w} q(\frac {\w} {N}) h(y,\w)\right 
)^{\wedge}(\nu)\\
 &=& \int \frac 1 {\w} q(\frac {\w} {N}) h(y,\w) dt +
 \sum\limits_{\nu \neq 0} \hat H _{1, N} (y, \nu) + \sum\limits_{\nu \neq 0} 
\hat H_{2, N} (y, \nu)
 \end{eqnarray*}
 where 
 $$H_{i, N} (y, t) = \frac 1 {\w} q(\frac {\w} {N}) h_i(y,\w), \;\;\;  i = 1,\; 
2.$$
Applying {\bf Lemma 1} and {\bf Lemma 2} with {\bf Remark 2} we bound the sum:
\begin{eqnarray*}
\|\sum\limits_{l 
\ge 0}^{k}\int\limits_{0}^{\infty} q(\frac {t} {2^l}) h(y,t) dt\|_{p'} 
&\le&
\sum\limits_{l \ge 0}\sum\limits_{\nu \neq 0} \|\hat H _{1, 2^l} (y, \nu)\|_{p'}
+ \sum\limits_{\nu \neq 0} 
\|\sum\limits_{l \ge 0}^{k}\hat H_{2, 2^l} (y, \nu)\|_{p'}\\
&\le&
C\|f\|_p.
\end{eqnarray*}
Combining \bref{convol}, \bref{zeroter} and the last inequality
we obtain the desired result
\begin{eqnarray*}
\|f*\widehat {\psi_k} 
\|_{p'} \le C\|f\|_p
\end{eqnarray*}
from which the statement of the theorem follows.
\hfill$\square$
\begin{rem}
We say that a function $f \in L^p$ has vanishing periodizations if there exists 
a sequence of Schwartz functions $f_k$ with vanishing periodizations converging 
to $f$ in $L^p$. It follows from {\bf Theorem 1} that $f \in L^{p'}$ and $f_k$ 
converge to $f$ in $L^ {p'}$ if dimension $d \ge 3$ and $1 \le p 
< \frac {2d}{d + 2}$.
\end{rem}

\section{Counterexamples and open questions}\label{counterex}
{\bf Theorem 1} does not say what happens when $d = 1$ and $d = 2$.\\

$d = 1$ is not an interesting case. We can easily construct 
examples of functions $f$ with vanishing periodizations such that their $L^p$ 
norms are not bounded by their $L^q$ norms for any given pair of $p \neq q$.\\

When $d = 2$ {\bf Theorem 1} does not hold. More precisely, {\bf Lemma 3} below 
shows that if $1 \le p < 2$ 
then
the following inequality does not hold for functions with vanishing 
periodizations:
$$\|f\|_{p'} \lesssim \|f\|_p.$$
In this lemma we  will deal with a sequence of
functions $f_n$ such that
$\hat f_n$
 vanish on all circles of radius $\sqrt{l^2 + k^2}$. Denote by $X_2$ the Banach 
space 
of functions from $L^1({\Bbb {R}} ^2)$ whose Fourier transforms vanish on all 
circles of radius $\sqrt{l^2 + k^2}$
$$X_2 = \{f \in L^1({\Bbb {R}} ^2): \hat f({\bf r}) = 0 \text { if } |{\bf r}| = 
\sqrt{l^2 + k^2},
 (k,l) 
\in {\Bbb {Z}} ^2 \}.$$ 
The next lemma crucially depends on the 
following fact from the Number Theory (\cite{G}, p.22):\\
{\it The number of integers in $[n, 2n]$ which can be represented as sums of two 
squares is
$n \epsilon_n$ where $\epsilon_n \lesssim \frac {1}{\ln^{1/2} n} 
\rightarrow 0$ as $n \rightarrow \infty$. }\\
We only use the fact that $\lim \epsilon_n = 0$.

\begin{lem} Let $1 \le p < 2$ and $d = 2$ then there exists a sequence of 
Schwartz functions $f_n \in X_2$ such that
\begin{eqnarray*}\lim\limits_{n \rightarrow \infty}\frac 
{\|f_n\|_{p'}}{\|f_n\|_p} = 
\infty.
\end{eqnarray*} 

\end{lem}
{\bf Proof of Lemma 3:} Let $a_1 < a_2 < a_3 <... $ be the enumeration of 
numbers 
$a_m = \sqrt{l^2 + k^2}$ in ascending order. Denote $\delta_m = a_{m+1} - a_m$. 
As we 
already said
the number of $a_m$ in $[\sqrt n, 2\sqrt{n}]$ is $~n\epsilon_n$. Let $a_{m_0}$ 
and $a_{m_1}$ be
correspondingly the smallest and the largest such $a_m$. Then 
$$\sum\limits_{m=m_0}^{m_1 - 1} \delta_m = a_{m_1} - a_m  \sim \sqrt n.$$
Let 
\begin{eqnarray}
\delta = \frac C {\sqrt n \epsilon_n}\label{deltadef}
\end{eqnarray}
 with small enough constant $C > 0$ 
so 
that if
$$M = \{m, m_0 \le m < m_1: \delta_m \ge \delta\}$$ then
$$\sqrt n \lesssim \sum\limits_{m \in M} \delta_m $$
since $m_1 - m_0 \sim n\epsilon_n.$  Choose coordinate axes $x$ and $y$.
We will construct $\hat f_n$ supported in $\bigcup\limits_{m \in M} R_m$ where 
$R_m$ is
a largest possible rectangle inscribed between circles of radius $a_m$ and 
$a_{m+1}$ with sides
parallel to the coordinate axes. Then
$R_m$ is of size $ \sim \delta_m \times \sqrt{\delta_m a_m} \gtrsim \delta_m 
\times 
\sqrt{\delta \sqrt n}
\gtrsim \delta_m \times 
1.$  We will split each rectangle $R_m$ further into smaller 
$\left[\frac {\delta_m}{\delta}\right]$ rectangles $r$ of the same size 
$\sim \delta \times 1.$ The number of these rectangles $r$ is
\begin{eqnarray}
N &=& \sum\limits_{m \in M} \left[\frac 
{\delta_m}{\delta}\right]\nonumber \\ 
&\sim& \sum\limits_{m \in M} \frac {\delta_m}{\delta}\nonumber \\
 &\sim& \frac {\sqrt n}
{\frac 1 {\sqrt n \epsilon_n}} = n\epsilon_n \label{Ndef}
\end{eqnarray}
since $\delta_m \ge \delta$ for $m \in M.$ Enumerate these rectangles $r_k$, $k 
= 
1,...,N$. Let
$r_k$ be centered at $(\lambda_k,0)$  It is clear that $|\lambda_k -\lambda_l| 
\ge \delta$ for
$k \neq l$. Let $\phi$ be a nonnegative Schwartz function on ${\Bbb {R}}$ 
supported in 
$[-\frac 12, \frac 12]$.  We have that $\check \phi(x) \ge C > 0$ when $x$ is 
small enough. Define $\hat f_n$ as the following sum:
\begin{eqnarray}\hat f_n (x,y)= \sum\limits_{k=1}^{N}\phi(\frac {x - 
\lambda_k}{\delta})
\phi({y})\label{hatf_nd=2}.
\end{eqnarray}
The $k$-th term in \bref{hatf_nd=2} is supported in $r_k$. Therefore, $\hat f_n$ 
is a 
Schwartz function supported in $\bigcup\limits_{m \in M} R_m$.
Hence $\hat f_n$ vanishes on all circles of radius $a_l$. Taking the inverse 
Fourier transform of
 \bref{hatf_nd=2}, we get
\begin{eqnarray}f_n (\xi,\eta)= \delta\check\phi(\xi{\delta})
\check\phi(\eta )
\sum\limits_{k=1}^{N}e^{i\lambda_k \xi}.\label{f_nd=2}
\end{eqnarray}
Assume that $p' < \infty$.
Then \begin{eqnarray*} \int |f_n(\xi,\eta)|^{p'}
d\xi d\eta &\ge&
\|\check \phi \|_{p'}^{p'}\delta ^{p'}\int\limits_{|\xi| \le \frac 
{100^{-1}}{\sqrt {n}}}
|\check\phi(\xi{\delta})|^{p'} |\sum\limits_{k=1}^{N}e^{i\lambda_k 
\xi}|^{p'}d\xi\\
&\gtrsim& \delta ^{p'} N^{p'}\frac 1 {\sqrt {n}}\\
&\sim& (\sqrt {n})^{p' - 1}.
\end{eqnarray*}
To obtain the second inequality we used that
$$|\sum\limits_{k=1}^{N}e^{i\lambda_k 
\xi}| \ge |\sum\limits_{k=1}^{N} \cos{ (\lambda_k 
\xi)}| \gtrsim N$$ since $|\lambda_k 
\xi| \le \frac 1 {50}$.
We used \bref{deltadef} and \bref{Ndef} to obtain the last estimate.
Therefore
\begin{eqnarray}
\|f_n\|_{p'} \gtrsim (\sqrt {n})^{\frac 1 p} \label{f_nd=2infty}.
\end{eqnarray}
If $p' = \infty$ we can obtain in a similar way that
\begin{eqnarray}
\|f_n\|_{\infty} \ge |f_n(0)| \gtrsim \sqrt {n}. \label{p=infinity}
\end{eqnarray}
Now we will estimate the $L^p$ norm from above.
 Denote 
$$g(x) = \sum\limits_{k=1}^{N}e^{i\frac{\lambda_k}
{\delta} \xi}.$$
Since $|\frac {\lambda_k -\lambda_l}{\delta}| \ge \frac{\delta}{\delta} = 1$ for 
$k \neq l$ we
have 
\begin{eqnarray*}\int_{I}|g|^2 \sim N\end{eqnarray*}
for any interval $I$ of length $4\pi$ (see (\cite{Z}, Theorem 9.1)). Therefore,
\begin{eqnarray}\int_{I}|g|^p &\le& |I|^{1 - \frac 2 p}(\int_{I}|g|^2)^{\frac p 
2}\nonumber 
\\
&\lesssim&  N^{\frac p 2}\label{gd=2}\end{eqnarray}
for any interval $I$ of length $4\pi$.
Since $\check \phi$ is a Schwartz function, we have that
$$|\check \phi(x)| \lesssim \frac 1 {1 + x^2}.$$
Therefore
\begin{eqnarray}\int|f_n(\xi, \eta)|^p d\xi d\eta &=& 
\|\check \phi\|_p^{p} \delta ^{p-1}\int|\check \phi(\xi)|^p \cdot 
|\sum\limits_{k=1}^{N}e^{i\frac{\lambda_k}
{\delta} \xi}|^ p d \xi\nonumber \\
&=& C\delta ^{p-1}\sum\limits_{l=-\infty}^{\infty}\int\limits_{l4\pi}^ 
{(l+1)4\pi}
|\check \phi(\xi)|^p \cdot |g(\xi)|^p d \xi\nonumber \\
&\lesssim& \delta ^{p-1}\sum\limits_{l=-\infty}^{\infty}\frac 1 {(1 + l^2)^p} 
N^{\frac p 2}\nonumber 
\\
&\lesssim& \sqrt {n} \epsilon _n ^{1 - \frac p 2}\nonumber
.
\end{eqnarray}
We used \bref{deltadef} and \bref{Ndef} to obtain the last estimate. Therefore
\begin{eqnarray}
\|f_n\|_{p} \lesssim (\sqrt {n})^{\frac 1 p} \epsilon _n ^{\frac {2 -  p} {2p}}
\label{f_nd=2L1}
\end{eqnarray}
Dividing \bref{f_nd=2infty} by \bref{f_nd=2L1} we obtain the desired result
\begin{eqnarray*}\frac {\|f_n\|_{p'}}{\|f_n\|_p} &\ge& \frac {(\sqrt{n})^{\frac 
1 p}}
{(\sqrt {n})^{\frac 
1 p} \epsilon _n ^{\frac {2 -  p} {2p}}}
\\
&=& \frac 1 {\epsilon _n ^{\frac {2 -  p} {2p}}} \rightarrow \infty
\end{eqnarray*}
as $n \rightarrow \infty$
 since $p < 2$.
\hfill$\square$

\begin{cor}
There exists a function $f \in X_2$ such that 
$$\|f\|_{L^{\infty}(D(0,1))} = 
\infty.$$ 
\end{cor}
It follows immediately from the lemma and \bref{p=infinity} that
if $p = 1$ then
\begin{eqnarray*}\sup\limits_{f \in X_2}\frac 
{\|f\|_{L^{\infty}(D(0,1))}}{\|f\|_1}
= \infty.\end{eqnarray*}
We claim that there exists a function $f \in X_2$ such that 
$\|f\|_{L^{\infty}(D(0,1))} = 
\infty.$ Suppose towards a contradiction that this is not true. Then the 
restriction 
operator 
$$T: f \rightarrow f|_{D(0,1)}$$ maps $X_2$ to $L^{\infty}(D(0,1))$. 
Note that if 
$f_n \rightarrow f$ in $L^1$ and $f_n \rightarrow g$ in 
$L^{\infty}(D(0,1))$, then $f = g$ a.e. on $D(0,1)$. An application of the 
Closed 
Graph Theorem
shows that $T$ is a bounded operator acting from $X_2$ to $L^{\infty}(D(0,1)).$ 
This contradicts
to the {\bf Corollary 2}. Thus we proved our claim. 
\hfill$\square$\\

Obviously, this function $f$ 
is 
not continuous. Therefore, it can serve as a counterexample to the theorem of M. 
Kolountzakis and T. Wolff (\cite{KW}, Theorem 1) mentioned in {\bf Introduction} 
when $d = 2$.
\begin{rem}
However, it is not known whether the following inequality holds for $f \in X_2$:
$$\|f\|_r \lesssim \|f\|_p$$
where $1 \le p < 2$ and $p < r < p'$.
\end{rem}

Now we will show that the range of $r$ in {\bf Corollary 1} is sharp. We need to 
check two cases: $r > p'$ and $r < p$. In the former case the argument will be 
similar to the one in the previous lemma. Therefore we will give only a sketch 
of the proof. We  will deal with a sequence of
functions $f_n$ such that
$\hat f_n$
 vanish on all circles of radius $\sqrt{m_1^2 + ... + m_d^2}$. Denote by $X_d$ 
the 
Banach 
space 
of functions from $L^1(\r)$ whose Fourier transforms vanish on all 
circles of radius $\sqrt{m_1^2 + ... + m_d^2}$
$$X_d = \{f \in L^1(\r): \hat f({\bf r}) = 0 \text { if } |{\bf r}| = 
\sqrt{m_1^2 + ... + m_d^2},
 (m_1, ..., m_d) 
\in \z \}.$$ 
We will construct a sequence of Schwartz functions $f_n$ with Fourier transforms  
supported outside of spheres of radius $\sqrt m$. Therefore these functions 
automatically belong to $X_d$.
\begin{lem} Let $1 < p \le 2$ and $r > p'$ then there exists a sequence of 
Schwartz functions $f_n \in X$ such that
\begin{eqnarray*}\lim\limits_{n \rightarrow \infty}\frac 
{\|f_n\|_r}{\|f_n\|_{p}} = 
\infty.
\end{eqnarray*} 

\end{lem}
{\bf Proof of Lemma 4:} A maximal rectangle inscribed between spheres of radius 
$\sqrt n$ and $\sqrt {n + 1}$  has dimensions $\sim \frac 1 {\sqrt n} \times 1 
\times 1 \times ... \times 1$. Let $r_k$ be parallel identical rectangles inscribed 
between spheres of radius $\sqrt {n + k}$ and $\sqrt {n + k + 1}$, where $k = 0, 
1, ..., n - 1$, with dimensions $\sim \frac 1 {\sqrt n} \times 1 \times 1 \times 
... \times 1$ and centered at $(\lambda_k, 0, 0, ..., 0)$. It is clear that 
$\lambda_{k + 1} - \lambda_k \sim \frac 1 {\sqrt n}$. Let $\phi$ be a 
nonnegative Schwartz function on ${\Bbb {R}}$ 
supported in 
$[-\frac 1 {100}, \frac 1 {100}]$.  We have that $\check \phi(x) \ge C > 0$ when 
$x$ is 
small enough. Define $\hat f_n$ as the following sum:
\begin{eqnarray}\hat f_n (x_1,x_2, ..., x_d)= \sum\limits_{k=0}^{n-1}\phi( (x_1 
- 
\lambda_k)\sqrt n)
\prod\limits_{l = 2}^{d}\phi(x_l)\label{hatf_nd=d}.
\end{eqnarray}
The $k$-th term in \bref{hatf_nd=d} is supported in $r_k$. Therefore, $\hat f_n$ 
is a 
Schwartz function vanishing on all spheres of radius $\sqrt m$. Taking the 
inverse 
Fourier transform of
 \bref{hatf_nd=d}, we get
\begin{eqnarray}f_n (y_1,y_2,..., y_d)= \prod\limits_{l = 2}^{d}\check\phi(y_l)
\frac 1 {\sqrt n}\check\phi(\frac {y_1} {\sqrt 
n})\sum\limits_{k=0}^{n-1}e^{i\lambda_k y_1}.\label{f_nd=d}
\end{eqnarray}
Arguments analogous to those in {\bf Lemma 3} show that
\begin{eqnarray*}
\|f_n\|_r \gtrsim (\sqrt n)^{\frac 1 {r'}}
\end{eqnarray*}
and
\begin{eqnarray*}
\|f_n\|_p \lesssim (\sqrt n)^{\frac 1 p}.
\end{eqnarray*}
Therefore
\begin{eqnarray*}
\frac {\|f_n\|_r}{\|f_n\|_p} \gtrsim (\sqrt n)^{\frac 1 {p'} - \frac 1 {r}} 
\rightarrow \infty
\end{eqnarray*}
as $n \rightarrow \infty$ since $r > p'$.
\hfill$\square$\\

The case when $r < p$ is very simple. Let 
$$\hat f (x) = \phi (\frac {x - x_0}{\epsilon})$$
where $\phi$ is a Schwartz function supported in $B^d(0, 1)$ so that $\hat f$ is 
supported in a small ball $B^d(x_0, \epsilon)$ placed between two fixed spheres 
of 
radius $\sqrt n$ and $\sqrt {n + 1}$. Then $f (y) = \epsilon ^d \check \phi 
(\epsilon y)$ and
$$\frac {\|f\|_r} {\|f\|_p} \sim \frac {\epsilon ^{\frac d {r'}}} {\epsilon 
^{\frac d {p'}}} \rightarrow \infty$$
as $\epsilon \rightarrow 0$ since $r < p$. Note that we didn't put any 
restriction on $p$ here.\\

Now we will show that {\bf Theorem 1} does not hold if $p > 2$. More precisely,
let $p > 2$ and $ r \neq p$ then the following inequality is not true for 
functions with vanishing periodizations:
$$\|f\|_r \lesssim \|f\|_p.$$
We just considered the case when $r < p$ therefore we need to consider only the 
case $r > p$. The argument is almost the same as in {\bf Lemma 4}. We can 
construct a sequence of Schwartz functions $f_n$ with Fourier transforms 
vanishing on all spheres of radius $\sqrt m$ and such that $\|f_n\|_r \gtrsim 
(\sqrt n)^{\frac 1 {r'}}$ and $\|f_n\|_p \le \|\hat f_n\|_{p'} \lesssim (\sqrt 
n)^{\frac 1 {p'}}$. Therefore
$$\frac {\|f_n\|_r} {\|f_n\|_p} \gtrsim (\sqrt n)^ {\frac 1 p - \frac 1 r} 
\rightarrow \infty.$$

\begin{rem}
Since {\bf Theorem 1} trivially holds for $p = 2$  it is natural to expect that
it should hold for $1 \le p \le 2$. It is unknown whether the {\bf Theorem 1} 
holds for $\frac {2d}{d + 2} \le p < 2$.
\end{rem}

Another interesting question is whether the following is true:
\begin{eqnarray}
\|\hat f\|_p \lesssim \|f\|_p \label{Fourierconj}
\end{eqnarray}
for some range of $p < 2$ if $f$ has vanishing periodizations. It would then 
follow that
\begin{eqnarray}
\|\hat f\|_r \lesssim \|f\|_p \label{Fourier}
\end{eqnarray}
for $p \le r \le p'$.
All we know from {\bf Theorem 1} is that
\bref{Fourier} holds
when $2 \le r \le p'$, $1 \le p < \frac {2d}{d + 2}$ and $d \ge 3$ since 
$\|f\|_2 \lesssim \|f\|_p$.\\

Our final open question is whether the following inequalities are true for 
functions with not necessarily vanishing periodizations $g_{\rho}$:

\begin{eqnarray*}
\|f\|_{p'} \lesssim \|f\|_p + \|g\|_{p'}
\end{eqnarray*}
and
\begin{eqnarray*}
\|g\|_{p'} \lesssim \|f\|_p + \|f\|_{p'}
\end{eqnarray*}
for some range of $p \le \frac {2d}{d + 1}$
where
\begin{eqnarray*}
\|g\|_{p'} = \left ( \int\limits_{\rho \in SO(d)} \|g_{\rho}\|_{p'}^p 
d\rho\right)^{\frac 1 p}.
\end{eqnarray*}

~

\end{document}